\documentclass[a4paper,11pt]{article}
\usepackage{tipa}
\usepackage{bbm}
\usepackage{bbding}

\usepackage{stmaryrd}
\usepackage{amsfonts}
\usepackage{graphicx}
\usepackage{epsfig}
\usepackage{amssymb}
\usepackage{amssymb,amsmath,amsthm}
\usepackage{pifont}
\usepackage{graphicx}
\usepackage{amsmath}
\usepackage{amssymb,amsmath,amsthm}
\usepackage{indentfirst}
\usepackage{cases}
\large\normalsize
\def\dis{\displaystyle}

\numberwithin{equation} {section} \makeatletter
\renewcommand{\@seccntformat}[1]{\csname
the#1\endcsname.\hspace{0.5em}} \makeatother

\begin{document}

 \title{\bf Vanishing Pressure Limit of Solutions to the Aw-Rascle
 Model for Modified Chaplygin gas
 \thanks{\ Supported by the NSF of China (11361073).} }

\author{Jinhuan Wang \ \ \  Jinjing Liu \ \ \  Hanchun Yang
\\ \footnotesize\slshape{Department of Mathematics, Yunnan University, Kunming, Yunnan 650091, P.R.China}
 }

\date{}\maketitle

{\small

\noindent{\bf Abstract:} This paper analyzes the vanishing pressure limit of
solutions to the Aw-Rascle model and the perturbed Aw-Rascle model for modified Chaplygin gas.
Firstly, the Riemann problem of the Aw-Rascle model is solved constructively.
A special delta shock wave in the limit of Riemann solutions is obtained.
Secondly, the Riemann problem of the perturbed Aw-Rascle model is solved
analytically. It is proved that, as the pressure tends to zero,
any Riemann solution containing two shock wave tends to a delta
shock solution to the transport equations; any Riemann solution
containing two rarefaction wave tends to a two-contact-discontinuity
solution to the transport equations and the nonvacuum intermediate
state in between tends to a vacuum state.

\vspace{0.1cm}

}

\date{}

\section{ Introduction }

The Aw-Rascle model of traffic flow reads
\begin{align}\label{eq1.1}
\left\{\begin{array}{l}
 \rho_t+(\rho u)_x=0,\cr\noalign{\vskip2truemm}
 \big(\rho (u+P)\big)_t+\big(\rho u(u+P)\big)_x=0,
\end{array}\right.
\end{align}
where $\rho\geq0$ and $u\geq0$ represent the traffic density and
velocity, respectively.  $P$ is the velocity offset and called as
the ``pressure" inspired from gas dynamics. The system \eqref{eq1.1} was proposed
by Aw and Rascle \cite{A-R} to remedy the deficiencies of second
order models of car traffic pointed out by Daganzo \cite{Daganzo}
and had also been independently derived by Zhang \cite{Zhang}. As a
macroscopic system, \eqref{eq1.1} is widely used to study the formation and
dynamics of traffic jams. The Riemann solutions of \eqref{eq1.1} with the
classical pressure
\begin{equation}\label{eq1.2}
P(\rho)=\rho^{\gamma},\gamma>0
\end{equation}
were obtained at low densities in \cite{A-R}. Lebacque, Mammar
and Salem \cite{Aw-Klar} solved the Riemann problem of \eqref{eq1.1} and
\eqref{eq1.2} with  extended fundamental diagram (equilibrium speed-density
relationship) for all possible initial data. Shen and
Sun \cite{Shen-Sun} considered \eqref{eq1.1} with a perturbed pressure term
\begin{equation}\label{eq1.3}
P(\rho)=\varepsilon\rho^{\gamma}, \varepsilon>0.
\end{equation}
The formal limit of \eqref{eq1.1} with \eqref{eq1.3} when
$\varepsilon\rightarrow 0$  is the transport equations
\begin{align}\label{eq1.4}
\left\{\begin{array}{l}
 \rho_t+(\rho u)_x=0,\cr\noalign{\vskip2truemm}
 (\rho u)_t+(\rho u^2)_x=0,
\end{array}\right.
\end{align}
also called the one-dimensional system of pressureless Euler
equations which has been studied and analyzed extensively since
1994. With radon measure as initial data, Bouchut \cite{B}
first presented an explicit formula of the Riemann solution and
checked the solution satisfying \eqref{eq1.4} in the sense of measure.
E, Rykov and Sinai \cite{WE} discussed the behaviour of global
weak solutions with random initial data. Sheng and Zhang
\cite{Sheng-Zhang} completely solved the 1-D and 2-D Riemann
problems of \eqref{eq1.4}. It has been shown that delta shock waves and vacuum states appear
in Riemann solutions of \eqref{eq1.4}. As for delta shock waves, we refer to
\cite{Tan-Zhang-Zheng,Li-Zhang,Yang,Yang-Zhang1,Yang-Zhang2}, etc.

For the researches of the delta shocks, one very interesting topic is to analyze the
formation of delta shock waves and vacuum states in
solutions. In \cite{C-L-1}, Chen and Liu considered the limit behaviour of Riemann solutions to
the Euler equations of isentropic gas dynamics for polytropic gas pressure $p(\rho)=\rho^{\gamma}/\gamma$ with $\gamma>1$ as the pressure vanishes. Further, they generalized this
result to the nonisentropic fluids in \cite{C-L-2}. Specially, Li \cite{Li} studied the pressure vanishing limit of solutions to the isentropic Euler equations for polytropic gas pressure when $\gamma=1$. For the related work, readers can see \cite{Sheng-Yin} for the relativistic Euler equations for polytropic gas by
Yin and Sheng, \cite{Yang-Wang1,Yang-Wang2} for the isentropic Euler equations for modified Chaplygin gas by Yang and Wang, and \cite{Cheng-Yang} for the Aw-Rascle model as the modified Chaplygin gas pressure tends to the Chaplygin gas pressure by Cheng and Yang, etc.

In \cite{Shen-Sun}, the authors  analyzed the limits of Riemann
solutions of \eqref{eq1.1} with \eqref{eq1.3} when $\varepsilon\rightarrow 0$. From
their results, one can see that the delta shock wave in the limit of Riemann
solutions dose not converge to that of \eqref{eq1.4}.  Recently, Pan and Han
\cite{Pan-Han} took the pressure
\begin{equation}\label{eq1.5}
P(\rho)=-\frac{\varepsilon}{\rho}
\end{equation}
for a Chaplygin gas in \eqref{eq1.1}.  They proved that the limits of Riemann
solutions of \eqref{eq1.1} with \eqref{eq1.5} are those of \eqref{eq1.4} when
$\varepsilon\rightarrow 0$.

In the present paper, we pay attention to the following pressure for the
modified Chaplygin gas
\begin{equation}\label{eq1.6}
P(\rho)=A\rho-\frac{B}{\rho^{\alpha}},\ (0<\alpha\leq1),
\end{equation}
where parameters $A,B$ and $\alpha$ are positive constants. Clearly,
when $A=0$, \eqref{eq1.6} corresponds to the generalized Chaplygin gas
\cite{Kamenshchik,Bento}, and if in addition $\alpha=1$, it becomes
the Chaplygin gas which was introduced by Chaplygin \cite{Ch} as a
suitable mathematical approximation for calculating the lifting
force on a wing of an airplane in aerodynamics. For Chaplygin gas,
we refer to \cite{Brenier,Guo,C-H-2,C-H-3}. If $B=0$, \eqref{eq1.6} is just the standard equation
of state of perfect fluid. The modified Chaplygin gas interpolates
between  Chaplygin gas fluids at low energy densities and standard
fluids at high energy densities. It was introduced by Benaoum
\cite{Benaoum-1,Benaoum-2} to describe the current accelerated
expansion of the universe at large cosmological scales and to use it
as a suitable kind of candidates of dark energy. Thus \eqref{eq1.1} with
\eqref{eq1.6} may model the motion process of the microscopic particles in
the universe. Obviously, when ${A,B\rightarrow0}$, \eqref{eq1.1} with \eqref{eq1.6} also
formally become the transport equations.

In this paper, we are concerned with limits of Riemann  solutions
of \eqref{eq1.1} with \eqref{eq1.6} as the parameters ${A,B\rightarrow0}$.
We prove that a special delta shock wave satisfying a special $\delta$-entropy condition develops
in the limit of Riemann solutions of \eqref{eq1.1} with \eqref{eq1.6} as the parameters ${A,B\rightarrow0}$. This delta shock wave is not exactly that of \eqref{eq1.4}.

In order to solve
it, we introduce a perturbations to model \eqref{eq1.1} with \eqref{eq1.6}
such that it becomes the following perturbed Aw-Rascle  model
\begin{align}\label{eq1.7}
\left\{\begin{array}{l}
 \rho_t+(\rho u)_x=0,\cr\noalign{\vskip2truemm}
 \Bigg(\rho \Big(u+\dis\frac{A}{2}\rho-\frac{1}{1-\alpha}\frac{B}{\rho^{\alpha}}\Big)\Bigg)_t+\Bigg(\rho u\Big(u+A\rho-\frac{B}{\rho^{\alpha}}\Big)\Bigg)_x=0,
\end{array}\right.\quad
0<\alpha<1.
\end{align}
It is proved that the delta shock wave and vacuum state appear in the
limits of Riemann solutions of \eqref{eq1.7} when ${A,B\rightarrow0}$, which are exactly the solutions to the transport equations. It is also noticed
that \eqref{eq1.7} fails as $\alpha=1$. For this
case, we will consider it in the future.

Firstly, we solve the Riemann problem of \eqref{eq1.1} and \eqref{eq1.6} with
initial data
\begin{align}\label{eq1.8}
(u,\rho)(0,x)=\left\{\begin{array}{ll}
  (u_-,\rho_-),&x<0,\cr\noalign{\vskip1truemm}
  (u_+,\rho_+), &x>0,
 \end{array}\right.
\end{align}
where $(u_\pm,\rho_\pm)$ are arbitrary constants and $u_\pm>0,
\rho_\pm>0$. Since one eigenvalue is genuinely nonlinear and the
other is linearly degenerate, the elementary waves consist of
rarefaction wave ($R$), shock wave ($S$) and contact discontinuity
($J$). The curves of elementary waves divide the phase plane into
four domains. By the phase plane analysis method, we establish
the existence and uniqueness of Riemann solutions including two
different structures $R+J$ and $S+J$. Then, it is shown that, as $A,B\rightarrow0$,
when $u_+<
u_-$, the Riemann solution $S+J$ converges to a special delta shock
solution, whose propagation speed and strength are different from
that of the transport equations. Besides, it is also shown that when
$u_+> u_-$, the Riemann solution $R+J$ tends to a
two-contact-discontinuity solution to transport equations, and the
nonvacuum intermediate state between $R$ and $J$ tends to a vacuum
state.

Secondly,  we solve the Riemann problem \eqref{eq1.7} and \eqref{eq1.8}. The
elementary waves contain the backward (forward) rarefaction wave
$\overleftarrow{R}$ ($\overrightarrow{R}$) and  backward (forward)
shock wave $\overleftarrow{S}$ ($\overrightarrow{S}$). These curves
of elementary waves divide the phase plane into four regions. With
the phase plane analysis, we obtain four kinds of Riemann solutions:
$\overleftarrow{R}\overrightarrow{R},
\overleftarrow{R}\overrightarrow{S},
\overleftarrow{S}\overrightarrow{R}$,$\overleftarrow{S}\overrightarrow{S}$.
Then, we prove that when $u_+< u_-$ and $A,B\rightarrow0$, the
Riemann solution containing two shock waves exactly converges to a delta
shock solution to the transport equations. We also
prove that, when $u_+> u_-$ and $A,B\rightarrow0$, the Riemann
solution containing two rarefaction waves tends to a
two-contact-discontinuity solution to the transport equations, and the
nonvacuum intermediate state between the two rarefaction waves tends
to a vacuum state.

The organization of this paper is as follows. In Section 2, we
review the Riemann solutions of \eqref{eq1.4} and \eqref{eq1.8}.  Section 3 solves the
Riemann problem \eqref{eq1.1}, \eqref{eq1.6} and \eqref{eq1.8}. Sections 4 and  5
investigate the limits of solutions of \eqref{eq1.1}, \eqref{eq1.6} and \eqref{eq1.8}. In
Section 6, we solve the Riemann problem \eqref{eq1.7} and \eqref{eq1.8}. Sections 7
and 8 analyze the limits of solutions of \eqref{eq1.7} and \eqref{eq1.8}.

 \section{Delta-shocks and vacuums for the transport equations}

For completeness, this section briefly recalls delta shock waves and
vacuum states in the Riemann solutions to the transport equations
\eqref{eq1.4}, see \cite{Sheng-Zhang} for more details.

System \eqref{eq1.4} has a double eigenvalue $\lambda=u$ with the associated
eigenvector $r=(0,1)^T$ satisfying $\nabla\lambda\cdot r\equiv0$,
which means the system \eqref{eq1.4} is nonstrictly hyperbolic and $\lambda$
linearly degenerate.

Consider the Riemann problem \eqref{eq1.4} and \eqref{eq1.8}. By seeking
self-similar solution $(u,\rho)(t,x)=(u,\rho)(\xi)\ (\xi=x/t)$, it
is easy to find that, besides the constant state and singular
solution $u=\xi,\rho=0$ (called vacuum states), the elementary waves
of \eqref{eq1.4} are nothing but the contact discontinuity. The Riemann
problem \eqref{eq1.4} and \eqref{eq1.8} can be solved by two cases.

For the case $u_- < u_+$, the solution includes two contact
discontinuities and a vacuum state besides constant states. That is,
\begin{align}\label{eq2.1}
(u,\rho)(\xi)=\left\{\begin{array}{lc}
 (u_-,\rho_-),&-\infty<\xi\leq u_-, \cr\noalign {\vskip1truemm}
 (\xi,0),&u_-<\xi<u_+, \cr\noalign {\vskip1truemm}
 (u_+,\rho_+),&u_+\leq\xi<+\infty.
 \end{array}\right.
\end{align}

For the case $u_- > u_+$, the solution is a delta-shock wave type one.

In order to define the measure solution, the weighted
$\delta$-function $w(s)$$\delta_S$ supported on a smooth curve $S$
parameterized as $t=t(s)$, $x=x(s) (c\leq s\leq d)$ can be defined
by
\begin{equation}\label{eq2.2}
\langle w(t(s))\delta_S, \varphi(t(s),x(s))\rangle =\dis\int^d_c\dis
w(t(s))\varphi(t(s),x(s))\sqrt{x'(s)^2 + t'(s)^2}ds
\end{equation}
for all  test functions $\varphi(t,x)\in C^\infty_0 (R^+ \times
R)$.

With this definition, a delta-shock solution of \eqref{eq1.4} can be expressed as
\begin{equation}\label{eq2.3}
\rho(t,x)=\rho_0(t,x)+w(t)\delta_S, \ \ u(t,x)=u_0(t,x),
\end{equation}
where $S=\{(t,\sigma t):0\leq t <\infty\}$,
\begin{equation}\label{eq2.4}
\rho_0(t,x)=\rho_-+[\rho]\chi(x-\sigma t),
u_0(t,x)=u_-+[u]\chi(x-\sigma t),
w(t)=\dis\frac{t}{\sqrt{1+\sigma^2}}(\sigma[\rho]-[\rho
u]),
\end{equation}
in which $[g]=g_+-g_-$ denotes the jump of function $g$ across the
discontinuity, $\sigma$ is the velocity of the delta-shock, and
$\chi(x)$ the characteristic function that is 0 when $x<0$ and 1
when $x>0$.

As shown in \cite{Sheng-Zhang}, for any $\varphi(t,x)\in C^\infty_0
(R^+ \times R)$, the delta-shock solution constructed above
satisfies
\begin{align}\label{eq2.5}
\begin{array}{l}
\langle \rho, \varphi_t\rangle +\langle \rho
u,\varphi_x\rangle=0,\cr\noalign{\vskip2truemm} \langle \rho
u,\varphi_t\rangle+\langle \rho u^2,
\varphi_x\rangle=0,\cr\noalign{\vskip2truemm}
\end{array}
\end{align}
where
\begin{align}\label{eq2.6}
\begin{array}{l}
\langle \rho,
\varphi\rangle=\dis\int^{+\infty}_0\dis\int^{+\infty}_{-\infty}\dis\rho_0
\varphi dx dt+\langle
w\delta_S,\varphi\rangle,\cr\noalign{\vskip2truemm} \langle\rho
u,\varphi\rangle=\dis\int^{+\infty}_0\dis\int^{+\infty}_{-\infty}\dis\rho_0
u_0 \varphi dx dt+\langle \sigma w\delta_S,\varphi\rangle.
\end{array}
\end{align}

Under the above definitions, the
generalized Rankine-Hugoniot relation reads
\begin{align}\label{eq2.7}
\left\{\begin{array}{l}
\dis\frac{dx}{dt}=\sigma,\cr\noalign{\vskip2truemm}
\dis\frac{d\big(w(t)\sqrt{1+\sigma^2}\big)}{dt}=\sigma[\rho]-[\rho
u],\cr\noalign{\vskip2truemm}
\dis\frac{d\big(w(t)\sigma\sqrt{1+\sigma^2}\big)}{dt}=\sigma[\rho
u]-[\rho u^2],
\end{array}\right.
\end{align}
which reflects the relationships among the location, weight and
propagation speed of the delta shock wave.

To guarantee the uniqueness, the entropy condition is supplemented
as
\begin{equation}\label{eq2.8}
u_+< \sigma< u_-,
\end{equation}
which means that all the characteristic lines on both sides of the
discontinuity are not out-going. So it is a overcompressive shock
wave.

Solving the equations \eqref{eq2.7} with
initial data $x(0)=0$ and $w(0)=0$ under the entropy condition \eqref{eq2.8}
yields
\begin{equation}\label{eq2.9}
\sigma=\dis\frac{\sqrt{\rho_+}u_+ +
\sqrt{\rho_-}u_-}{\sqrt{\rho_+}+\sqrt{\rho_-}}\ \ \   \text{and}\ \
\ w(t)=\dis\frac{\sqrt{\rho_+ \rho_-}(u_- - u_+)t}{\sqrt{1+\sigma
^2}}.
\end{equation}
Therefore, a delta shock solution defined by \eqref{eq2.3} with \eqref{eq2.4} and
\eqref{eq2.9} is obtained.

\section{Solutions of Riemann problem \eqref{eq1.1}, \eqref{eq1.6} and
\eqref{eq1.8}}

In this section, we solve the elementary waves and  construct the
solutions of Riemann problem \eqref{eq1.1}, \eqref{eq1.6} and
\eqref{eq1.8}.

For any
$A,B>0$, the system has two eigenvalues
\begin{align}\label{eq3.1}
\begin{array}{l}
\lambda_1=u-A\rho-\dfrac{B\alpha}{\rho^{\alpha}}\ ,\ \ \ \
\lambda_2=u
\end{array}
\end{align}
with right eigenvectors
$$
\begin{array}{l}
{r}_1=(-A-\dfrac{B\alpha}{\rho^{1+\alpha}}, 1)^T,\ \ \ \ {r}_2=(0, 1)^T
\end{array}
$$
satisfying
\begin{equation}\label{eq3.2}
 \nabla\lambda_1\cdot {r}_1=-2A-\dfrac{(1-\alpha)B\alpha}{\rho^{1+\alpha}}<0,\ \ \ \
 \nabla\lambda_2\cdot {r}_2\equiv0.
\end{equation}
Thus this system is strictly hyperbolic.  The first characteristic
is genuinely nonlinear and the associated wave is either shock wave
or rarefaction wave. The second is linearly
degenerate and  the  associated wave is the contact discontinuity.

Performing  the self-similar transformation $x/t=\xi$, we reach the
following boundary value problem
\begin{align}\label{eq3.3}
\left\{\begin{array}{l}
 -\xi\rho_\xi+(\rho u)_\xi=0,\cr\noalign{\vskip3truemm}
 -\xi\big(\rho (u+A\rho-\dfrac{B}{\rho^{\alpha}})\big)_\xi+\big(\rho u(u+A\rho-\dfrac{B}{\rho^{\alpha}})\big)_\xi=0,
\end{array}\right.
\end{align}
and
\begin{equation}\label{eq3.4}
(u,\rho)(\pm\infty)=(u_\pm,\rho_\pm).
\end{equation}

For any smooth solution, \eqref{eq3.3} is equivalent to
\begin{align}\label{eq3.5}
\left(\begin{array}{cc}
 \rho&u-\xi\cr\noalign {\vskip2truemm}
u-A\rho-\dfrac{B \alpha}{\rho^{\alpha}}-\xi&0
\end{array}\right)
\left(\begin{array}{c}
du\cr\noalign {\vskip2truemm}
  d\rho
\end{array}\right)=0,
\end{align}
which provides either the general solution (constant state)
\begin{equation}\label{eq3.6}
(u,\rho)(\xi)=constant,
\end{equation}
or rarefaction wave, which is a wave of the first characteristic
family,
\begin{align}\label{eq3.7}
R:\ \left\{\begin{array}{l}
 \xi=\lambda_1=u-A\rho-\dfrac{\alpha B}{\rho^{\alpha}},\cr\noalign {\vskip2truemm}
 u=-A\rho+\dfrac{B}{\rho^{\alpha}}+u_-+A\rho_--\dfrac{B}{\rho_-^{\alpha}},\ \ \ \rho<\rho_-.
\end{array}\right.
\end{align}

For a bounded discontinuity at  $\xi=\sigma^{AB}$,  the
Rankine-Hugoniot relation
\begin{align}\label{eq3.8}
\left\{\begin{array}{l}
 -\sigma^{AB}[\rho]+[\rho u]=0,\cr\noalign {\vskip3truemm}
 -\sigma^{AB}\Big[\rho (u+A\rho-\dfrac{B}{\rho^{\alpha}})\Big]+\Big[\rho
 u(u+A\rho-\dfrac{B}{\rho^{\alpha}})\Big]=0
\end{array}\right.
\end{align}
holds. By solving \eqref{eq3.8} and using the Lax
entropy inequalities, we obtain  the shock wave, which is a wave of
the first characteristic family,
\begin{align}\label{eq3.9}
S:\  \left\{\begin{array}{l}
 \sigma_1^{AB}=u-\dfrac{B}{\rho^{\alpha}}-A\rho_--B\dis\frac{(\rho_-^{1-\alpha}-\rho^{1-\alpha})}{\rho-\rho_-},\cr\noalign{\vskip3truemm}
u=-A\rho+\dfrac{B}{\rho^{\alpha}}+u_-+A\rho_--\dfrac{B}{\rho_-^{\alpha}},\
\ \ \ \ \ \ \rho>\rho_-.
\end{array}\right.
\end{align}

Besides, from \eqref{eq3.5} or \eqref{eq3.8}, one can easily get the contact
discontinuity, which is a wave of the second characteristic family,
\begin{equation}\label{eq3.10}
J:\
 \sigma_2^{AB}=u=u_-,\
\ \ \ \ \ \ \rho\lessgtr\rho_-.
\end{equation}

In the first quadrant of the $(u,\rho)$-plane,  the contact
discontinuity line $u=u_-$ is a straight one paralleling with the
$\rho$-axis. The rarefaction wave and the shock  wave curves have
the same expression
$u=-A\rho+\frac{B}{\rho^{\alpha}}+u_-+A\rho_--\frac{B}{\rho_-^{\alpha}}$,
which means the system belongs to the  Temple type \cite{Temple}.
Due to $u_\rho=-A-\frac{B\alpha}{\rho^{1+\alpha}}<0$ and
$u_{\rho\rho}=\frac{B\alpha(1+\alpha)}{\rho^{2+\alpha}}>0$, the two curves
are monotonic decreasing and convex. Moreover, it can be verified
that $\lim_{\rho\rightarrow0^+}u=+\infty$, which implies that the
rarefaction wave curve has the $u$-axis as the asymptote. It also
can be proved that $\lim_{\rho\rightarrow+\infty}u=-\infty$, which
implies the shock wave curve intersects with the $\rho$-axis at some
point. Fixing a left state $(u_-,\rho_-)$, the phase plane
 can be divided into four regions by the wave
curves, denoted by $I(u_-,\rho_-)$, $II(u_-,\rho_-)$,
$III(u_-,\rho_-)$ and $IV(u_-,\rho_-)$, respectively (see Fig. 1).

\begin{center}
\includegraphics*[160,496][351,632]{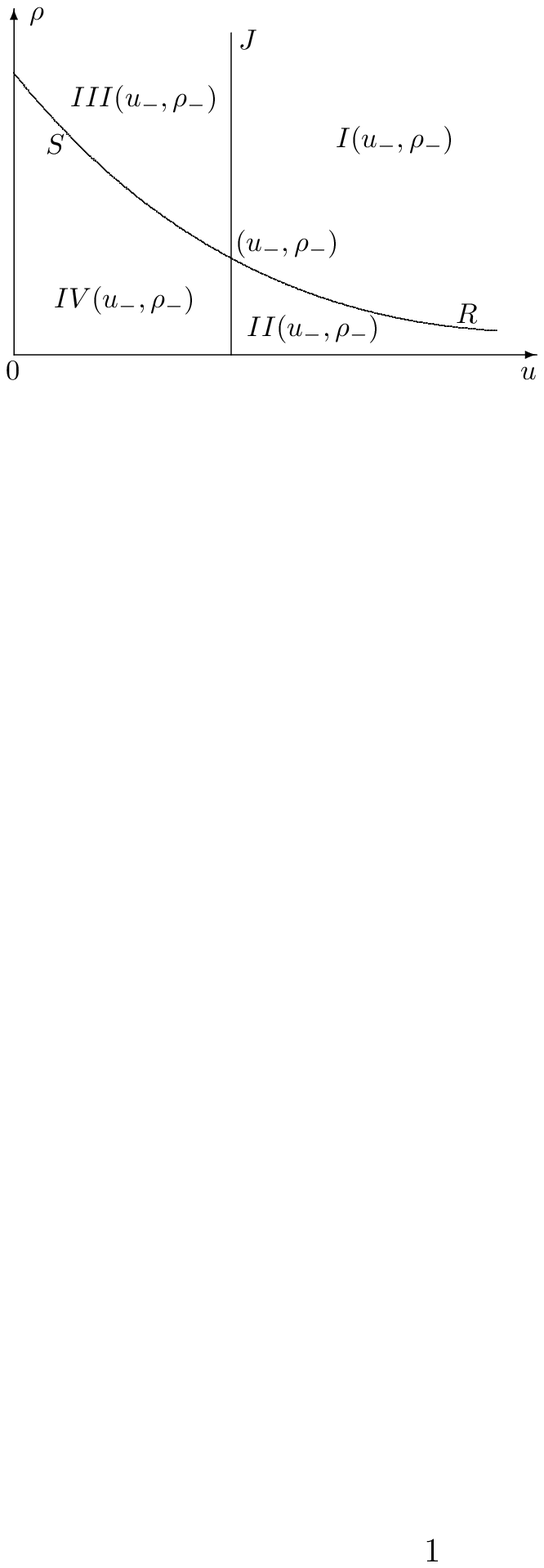}

{\text{ Fig.1.} Curves of elementary waves.}

\end{center}

Now, with phase plane analysis, according to the right
state $(u_+,\rho_+)$ in the different regions, one can get two kinds
of configurations of solutions:

$\textcircled{1}.\ (u_+,\rho_+)\in I(u_-,\rho_-)\cup
II(u_-,\rho_-):R+J$,

$\textcircled{2}.\ (u_+,\rho_+)\in III(u_-,\rho_-)\cup
IV(u_-,\rho_-):S+J$.

 \section{Formation of delta-shocks in solutions of \eqref{eq1.1}, \eqref{eq1.6} and \eqref{eq1.8}}

This section analyzes the limits as $A,B\rightarrow0$ of the Riemann
solutions of \eqref{eq1.1} with \eqref{eq1.6} and \eqref{eq1.8} when $u_+< u_-$. There are
two cases $(u_+,\rho_+)\in IV(u_-,\rho_-)$ and $(u_+,\rho_+)\in
III(u_-,\rho_-)$ to be considered.

\vspace{0.2cm}

We first discuss the case $(u_+,\rho_+)\in
IV(u_-,\rho_-)$.

\subsection{Limit behaviour of the Riemann solutions  as $A,B\rightarrow0$}

For any fixed $A,B>0$, when $u_+< u_-$ and  $(u_+,\rho_+)\in
IV(u_-,\rho_-)$, the solution of Riemann problem \eqref{eq1.1}, \eqref{eq1.6} and \eqref{eq1.8} is, besides two
constant states, a shock wave $S$ followed by a contact discontinuity $J$
with the intermediate state $(u_*^{AB},\rho_*^{AB})$. They have the
following relations
\begin{align}\label{eq4.1}
S:\  \left\{\begin{array}{l}
 \sigma_1^{AB}=u_*^{AB}-\dfrac{B}{(\rho_*^{AB})^{\alpha}}-A\rho_--B\dis\frac{(\rho_-^{1-\alpha}-(\rho_*^{AB})^{1-\alpha})}{\rho_*^{AB}-\rho_-},\cr\noalign{\vskip3truemm}
u_*^{AB}=-A\rho_*^{AB}+\dfrac{B}{(\rho_*^{AB})^{\alpha}}+u_-+A\rho_--\dfrac{B}{(\rho_-)^{\alpha}},\
\ \ \ \ \ \ \rho_*^{AB}>\rho_-,
\end{array}\right.
\end{align}
and
\begin{equation}\label{eq4.2}
 {J}:
 \sigma_2^{AB}=u_*^{AB}=u_+, \
\ \ \ \ \ \ \rho_*^{AB}>\rho_+,
\end{equation}
Then we have the following lemmas.

\vspace{0.2cm}

 \noindent{\textbf{Lemma 4.1.}}\textit{
$\lim\limits_{A,B\rightarrow0}\rho_*^{AB}=+\infty$.}

\vspace{0.2cm}

\textbf{Proof. } Suppose that
$\lim\limits_{A,B\rightarrow0}\rho_*^{AB}=a\in(\text{max}(\rho_-,\rho_+),+\infty)$.
It follows from \eqref{eq4.1} and \eqref{eq4.2} that
\begin{align}\label{eq4.3}
\begin{array}{l}
(\rho_*^{AB})^{\alpha}u_+=-A(\rho_*^{AB})^{1+\alpha}+B+(u_-+A\rho_--\dfrac{B}{\rho_-^{\alpha}})(\rho_*^{AB})^{\alpha},\
\ \
 \rho_*^{AB}>\rho_-.
\end{array}
\end{align}
Letting $A,B\rightarrow0$ in \eqref{eq4.3}, one can get $u_+=u_-$, which
contradicts  $u_+<u_-$. Therefore, Lemma 4.1 holds.

\vspace{0.2cm}

\noindent {\textbf{Lemma 4.2.}}\textit{
$\lim\limits_{A,B\rightarrow0}\sigma_1^{AB}=\lim\limits_{A,B\rightarrow0}\sigma_2^{AB}=u_+$.}

\vspace{0.2cm}

\textbf{Proof. } From \eqref{eq4.1}, \eqref{eq4.2} and Lemma 4.1, it is immediate
that
$$
\lim\limits_{A,B\rightarrow0}\sigma_1^{AB}=\lim\limits_{A,B\rightarrow0}u_*^{AB}=u_+.
$$

\vspace{0.2cm}

 Lemmas 4.1-4.2 show that when $A$ and $B$ drop to zero,  $S$ and $J$
coincide,  the intermediate density $\rho_*^{AB}$ becomes singular.

\vspace{0.2cm}

From \eqref{eq4.1}, \eqref{eq4.2} and Lemmas 4.1-4.2, we have

\vspace{0.2cm}

\noindent {\textbf{Lemma 4.3.}} \textit{$
\lim\limits_{A,B\rightarrow0}\rho_*^{AB}(\sigma_2^{AB}-\sigma_1^{AB})=\rho_-(u_--u_+).
$}

\subsection{Weighted delta shock waves }

Now, we give the theorem presenting the limit behaviour of Riemann solutions of \eqref{eq1.1}, \eqref{eq1.6} and \eqref{eq1.8} as
$A,B\rightarrow0$ for the case $u_+< u_-$ and $(u_+,\rho_+)\in
IV(u_-,\rho_-)$.

\vspace{0.2cm}

\noindent {\textbf{Theorem 4.4.}} \textit{Let  $u_+< u_-$ and
$(u_+,\rho_+) \in IV(u_-,\rho_-)$.  Assume $(u^{AB},\rho^{AB})(t,x)$
is the  Riemann  solution $S+J$ of \eqref{eq1.1}, \eqref{eq1.6} and \eqref{eq1.8}
constructed in Section 3. Then
$$
\lim\limits_{A,B\rightarrow0}u^{AB}(t,x)=\left\{\begin{array}{ll}
 u_-,&x<u_+ t,\\[1mm]
 u_+, &x=u_+ t,\\[1mm]
 u_+, &x>u_+ t,
 \end{array}\right.
 $$
$\rho^{AB}(t,x)$ converges in the sense of distributions, and the
limit  is  the sum of a step function and a $\delta$-function
supported on $x=u_+ t$ with the weight \
\textit{$\dis\frac{t}{\sqrt{1+u_+^2}}(\rho_-(u_--u_+))$.}}

\vspace{0.2cm}

The proof of Theorem 4.4 is similar to that of  Theorem 7.4 below,
here we omit it.

\vspace{0.2cm}

\noindent \textbf{Remark 4.5.} From  Theorem 4.4, it can be seen
that the propagation speed and strength of the delta shock wave in the
limit are different from that of the transport equations solved in
\cite{Sheng-Zhang}. For the limit solution of \eqref{eq1.1}, \eqref{eq1.6} and \eqref{eq1.8},
the characteristics on the right side of the delta shock wave
parallel to it while the characteristics on the left side of it come
into the delta shock wave. Thus, this limit of solution of \eqref{eq1.1}, \eqref{eq1.6} and \eqref{eq1.8}
 is not the entropy solution of the transport equations
 since \eqref{eq2.8} fails.

Especially, if we replace \eqref{eq2.8} by the special $\delta$-entropy
condition: $u_+=\sigma_0<u_-$. Then, when $u_+<u_-$ and
$(u_+,\rho_+) \in IV(u_-,\rho_-)$, the limit of the Riemann solution
of \eqref{eq1.1}, \eqref{eq1.6} and \eqref{eq1.8} is
$$
\left\{\begin{array}{ll} x(t)=\sigma_0 t=u_+t,\\
w_0(t)=\dis\frac{t}{\sqrt{1+u_+^2}}(\rho_-(u_--u_+)),
\end{array}\right.
$$
which satisfies the generalized Rankine-Hugoniot condition \eqref{eq2.7} and
is also the entropy solution to the transport equations under the
special $\delta$-entropy condition.

\vspace{0.2cm}

Then we are in the position to consider the case $(u_+,\rho_+) \in
III(u_-,\rho_-)$.

\vspace{0.2cm}

\noindent {\textbf{Lemma 4.6.}} \textit{Let  $u_+< u_-$, for any
fixed $(u_+,\rho_+)$ satisfying  $(u_+,\rho_+) \in III(u_-,\rho_-)$,
there exist $A_0,$ such that $(u_+,\rho_+) \in IV(u_-,\rho_-)$ when
$0<A<A_0$ and $0<B<A_0$.}

\vspace{0.2cm}

\textbf{Proof.} When $u_+< u_-$, for any fixed $(u_+,\rho_+)$,
$(u_+,\rho_+) \in III(u_-,\rho_-)$ is equivalent to
\begin{equation}\label{eq4.4}
u_++A\rho_+-\dis\frac{B}{\rho_+^{\alpha}}
> u_-+A\rho_--\dis\frac{B}{\rho_-^{\alpha}}, \quad \rho_+>\rho_-,
\end{equation}
and $(u_+,\rho_+) \in IV(u_-,\rho_-)$ is equivalent to
\begin{equation}\label{eq4.5}
u_++A\rho_+-\dis\frac{B}{\rho_+^{\alpha}} <
u_-+A\rho_--\dis\frac{B}{\rho_-^{\alpha}}.
\end{equation}
Combining \eqref{eq4.4} with \eqref{eq4.5}, we set $A=B=A_0$ and solve the equation
\begin{equation}\label{eq4.6}
u_++A_0\rho_+-\dis\frac{A_0}{\rho_+^{\alpha}} =
u_-+A_0\rho_--\dis\frac{A_0}{\rho_-^{\alpha}}
\end{equation}
to get
\begin{equation}\label{eq4.7}
A_0=\dis\frac{u_--u_+}{(\rho_+-\rho_+^{-\alpha})-(\rho_--\rho_-^{-\alpha})}.
\end{equation}

\noindent This completes the proof of Lemma 4.6.

\vspace{0.2cm}

From the analysis above, we can obtain the following theorem which
is an immediate consequence of Theorem 4.4 and Lemma 4.6.

\vspace{0.2cm}

\noindent {\textbf{Theorem 4.7.}} \textit{Let  $u_+< u_-$ and
$(u_+,\rho_+) \in III(u_-,\rho_-)$. Then there exist $A_0>0$, when
$A>A_0$ and $B>A_0$,  there is no delta shock wave in solutions.
When $0<A<A_0$ and $0<B<A_0$, $(u_+,\rho_+) \in IV(u_-,\rho_-)$,
then letting $A,B\rightarrow0$, the limit is a delta shock solution
which is the same as that  in Theorem 4.4. }

 \section{Formation of vacuums in solutions of \eqref{eq1.1}, \eqref{eq1.6} and \eqref{eq1.8}}

In this section, we turn to the limit as $A,B\rightarrow0$ of the
Riemann solutions of \eqref{eq1.1}, \eqref{eq1.6} and \eqref{eq1.8} when $u_+> u_-$. There
are also two cases: $(u_+,\rho_+)\in I(u_-,\rho_-)$ and
$(u_+,\rho_+)\in II(u_-,\rho_-)$.

\vspace{0.2cm}

First of all, we consider the case $(u_+,\rho_+)\in I(u_-,\rho_-)$.

\vspace{0.2cm}

For any fixed $A,B>0$, when $u_+> u_-$  and $(u_+,\rho_+)\in
I(u_-,\rho_-)$, the solution of Riemann problem \eqref{eq1.1}, \eqref{eq1.6} and \eqref{eq1.8}
 is a rarefaction wave $R $ followed by  a contact
discontinuity $J$ with the intermediate state
$(u_*^{AB},\rho_*^{AB})$, besides two constant states. Thus we have
\begin{equation}\label{eq5.1}
u_*^{AB}=-A\rho_*^{AB}+\dis\frac{B}{(\rho_*^{AB})^{\alpha}}+u_-+A\rho_--\frac{B}{\rho_-^{\alpha}},\
\ \ \rho_*^{AB}<\rho_-
\end{equation}
on $R$, and
\begin{equation}\label{eq5.2}
 \sigma_2^{AB}=u_*^{AB}=u_+,\ \ \rho_*^{AB}<\rho_+
\end{equation}
on $J$.

Then  $\rho_*^{AB}$ satisfies
\begin{equation}\label{eq5.3}
 u_+=-A\rho_*^{AB}+\frac{B}{(\rho_*^{AB})^{\alpha}}+u_-+A\rho_--\frac{B}{\rho_-^{\alpha}},\ \ \ \rho_*^{AB}<\rho_-.
\end{equation}

Now, we can conclude the following result.

 \vspace{0.2cm}

\noindent {\textbf{Theorem 5.1.}} \textit{When  $u_+>u_-$ and
$(u_+,\rho_+) \in I(u_-,\rho_-)$, the vacuum state occurs as
$A,B\rightarrow0$. And the rarefaction wave $R$ and contact
discontinuity $J$ become two contact discontinuities connecting the
constant states $(u_\pm,\rho_\pm)$ and the vacuum $(\rho=0)$. }

\vspace{0.2cm}

\textbf{Proof.}  Suppose that
$\lim\limits_{A,B\rightarrow0}\rho_*^{AB}=b\in
(0,\text{min}(\rho_-,\rho_+))$. Taking the limit $A,B\rightarrow0$
in \eqref{eq5.3}, we have $u_+=u_-$, which contradicts with
$u_+>u_-$. We thus obtain that $\lim\limits_{A,B\rightarrow
0}\rho_*^{AB}=0$.

Besides, we can find that $\lim\limits_{A,B\rightarrow
0}\lambda_1(u_-, \rho_-)=u_-$. It means that the rarefaction wave
$R$ tends to a contact discontinuity $\xi=x/t=u_-$ when
$A,B\rightarrow 0$.

In summary, in this case, the Riemann solution $R+J$ tends to a
two-contact-discontinuity solution to the transport equations \eqref{eq1.4}
with the same Riemann data \eqref{eq1.8}, and the nonvacuum state between
$R$ and $J$ becomes a vacuum state when $A,B\rightarrow 0$.

 \vspace{0.2cm}

We proceed to study the case $(u_+,\rho_+)\in II(u_-,\rho_-)$. Similarly to Lemma 4.6, we have

 \vspace{0.2cm}

\noindent {\textbf{Lemma 5.2.}} \textit{Let  $u_+> u_-$, for any
fixed $(u_+,\rho_+)$ satisfying  $(u_+,\rho_+) \in II(u_-,\rho_-)$,
there exist $\overline{A_0}$ such that $(u_+,\rho_+) \in
I(u_-,\rho_-)$ when $0<A<\overline{A_0}$ and $0<B<\overline{A_0}$.}

\vspace{0.2cm}

\textbf{Proof.} When $u_+> u_-$, for any fixed $(u_+,\rho_+)$,
$(u_+,\rho_+) \in II(u_-,\rho_-)$ is equivalent to
$$
u_++A\rho_+-\dis\frac{B}{\rho_+^{\alpha}} <
u_-+A\rho_--\dis\frac{B}{\rho_-^{\alpha}}, \quad \rho_+<\rho_-,
$$
and $(u_+,\rho_+) \in I(u_-,\rho_-)$ is equivalent to
$$
u_++A\rho_+-\dis\frac{B}{\rho_+^{\alpha}} >
u_-+A\rho_--\dis\frac{B}{\rho_-^{\alpha}}.
$$
Setting $A=B=\overline{A_0}$ and solving the same equation \eqref{eq4.6},
one can get $\overline{A_0}$ which has the same expression with
$A_0$ shown in \eqref{eq4.7}. This completes the proof of Lemma 5.2.

 \vspace{0.2cm}

Combining Theorem 5.1 with Lemma 5.2, we have the following result.

\vspace{0.2cm}

\noindent {\textbf{Theorem 5.3.}} \textit{Let  $u_+> u_-$ and
$(u_+,\rho_+) \in II(u_-,\rho_-)$. Then there exist
$\overline{A_0}$, when $A>\overline{A_0}$ and $B>\overline{A_0}$,
there is no vacuum in  solutions. When $0<A<\overline{A_0}$ and
$0<B<\overline{A_0}$, $(u_+,\rho_+) \in I(u_-,\rho_-)$, then taking
$A,B \rightarrow0$, the  vacuum state occurs which is the same as
that in Theorem 5.1.}

\vspace{0.2cm}

In Sections 6-8 below, we solve the Riemann problem of the perturbed
Aw-Rasde model \eqref{eq1.7} and study the limit behaviour of the Riemann
solutions of \eqref{eq1.7} as $A,B\rightarrow0$.

\section{Solutions of Riemann problem \eqref{eq1.7} and \eqref{eq1.8}}

In this section, we solve  the Riemann problem \eqref{eq1.7} and \eqref{eq1.8}. For
any fixed $A,B>0$, the system has two eigenvalues
\begin{align}\label{eq6.1}
\begin{array}{l}
\overline{\lambda_1}=u-\sqrt{u(A\rho+\dfrac{B\alpha}{\rho^{\alpha}})}\ ,\ \ \
\ \overline{\lambda_2}=u+\sqrt{u(A\rho+\dfrac{B\alpha}{\rho^{\alpha}})}
\end{array}
\end{align}
with right eigenvectors
$$
\begin{array}{l}
\overline{{r}_1}=\Bigg(
-\dis\frac{\sqrt{u(A\rho+\frac{B\alpha}{\rho^{\alpha}})}}{\rho},1\Bigg)^T,\ \
\ \
\overline{{r}_2}=\Bigg(\dis\frac{\sqrt{u(A\rho+\frac{B\alpha}{\rho^{\alpha}})}}{\rho},1\Bigg)^T
\end{array}
$$
satisfying $\nabla\overline{\lambda_i}\cdot \overline{{r}_i}\neq0\
(i=1,2)$ for $u>0$ and $\Big(3A\rho +\frac{B\alpha}{
\rho^{\alpha}}(2-\alpha)\Big)\sqrt{u}\neq(A\rho+\frac{B\alpha}{\rho^{\alpha}})^\frac{3}{2}$.
Thus this system is strictly hyperbolic and both the characteristic
fields are genuinely nonlinear in the phase plane
($u>0,\rho>0$).

Seeking the self-similar solution, we obtain the following boundary
value problem
\begin{align}\label{eq6.2}
\left\{\begin{array}{l}
 -\xi\rho_\xi+(\rho u)_\xi=0,\cr\noalign{\vskip3truemm}
 -\xi\Bigg(\rho
\Big(u+\dis\frac{A}{2}\rho-\frac{1}{1-\alpha}\frac{B}{\rho^{\alpha}}\Big)\Bigg)_\xi+\Bigg(\rho
u\Big(u+A\rho-\frac{B}{\rho^{\alpha}}\Big)\Bigg)_\xi=0,
\end{array}\right.
\end{align}
and \eqref{eq3.4}.

For any smooth solution, \eqref{eq6.2} is equivalent to
\begin{align}\label{eq6.3}
D\left(\begin{array}{c}
 du\cr\noalign {\vskip2truemm}
  d\rho
\end{array}\right)=0,
\end{align}
where
$$
D=\left(\begin{array}{cc}
 \rho&-\xi+u\cr\noalign {\vskip2truemm}
-\xi\rho+2\rho
 u+A\rho^{2}-\dfrac{B}{\rho^{1-\alpha}}&
 -\xi (u+ A\rho- \dfrac{B}{\rho^{\alpha}})+u^2+2A\rho
 u-\dfrac{B(1-\alpha)u}{\rho^{\alpha}}
 \end{array}\right),
$$
which provides either the constant state or the backward rarefaction
wave
\begin{align}\label{eq6.4}
\overleftarrow{R}(u_-,\rho_-):\ \left\{\begin{array}{l}
 \xi=\overline{\lambda_1}=u-\sqrt{u(A\rho+\dfrac{B\alpha}{\rho^{\alpha}})},\cr\noalign {\vskip2truemm}
 \sqrt{u}-\sqrt{u_-}=-\dis\frac{1}{2}\int^{\rho}_{\rho_-} \frac{\sqrt{As+\frac{B\alpha}{s^{\alpha}}}}{s} ds,
\end{array}\right.
\end{align}
or the forward  rarefaction wave
\begin{align}\label{eq6.5}
\overrightarrow{R}(u_-,\rho_-):\ \ \left\{\begin{array}{l}
\xi=\overline{\lambda_2}=u+\sqrt{u(A\rho+\dfrac{B\alpha}{\rho^{\alpha}})},\cr\noalign
{\vskip2truemm}
 \sqrt{u}-\sqrt{u_-}=\dis\frac{1}{2}\int^{\rho}_{\rho_-} \frac{\sqrt{As+\frac{B\alpha}{s^{\alpha}}}}{s} ds.
\end{array}\right.
\end{align}

For the backward  rarefaction wave, by differentiating $u$ with
respect to $\rho$ in the second equation of \eqref{eq6.4}, it follows that
$
u_{\rho}=\frac{-\sqrt{u(A\rho+\frac{B\alpha}{\rho^{\alpha}})}}{\rho}<0$. For
the forward rarefaction wave, it is easy to see that $
u_{\rho}=\frac{\sqrt{u(A\rho+\frac{B\alpha}{\rho^{\alpha}})}}{\rho}>0$.

Through differentiating $\xi$ with respect to $\rho$ and $u$ in the
first equation of \eqref{eq6.4} and noticing
$u_{\rho}=\frac{u_\xi}{\rho_\xi}$, we have
\begin{equation}\label{eq6.6}
1=\Bigg(1-\frac{\sqrt{A\rho+\frac{B\alpha}{\rho^{\alpha}}}}{2
u}+\frac{A\rho-\frac{B\alpha^2}{
\rho^{\alpha}}}{2(A\rho+\frac{B\alpha}{\rho^{\alpha}})}\Bigg)u_{\xi}.
\end{equation}
Since
\begin{equation}\label{eq6.7}
\frac{A\rho-\frac{B\alpha^2}{
\rho^{\alpha}}}{2(A\rho+\frac{B\alpha}{\rho^{\alpha}})}=\frac{1}{2}-\frac{\frac{B\alpha}
{\rho^{\alpha}}+\frac{B\alpha^2}{\rho^{\alpha}}}{2(A\rho+\frac{B\alpha}
{\rho^{\alpha}})} \quad \text{and} \quad 0<\frac{\frac{B\alpha}
{\rho^{\alpha}}+\frac{B\alpha^2}{\rho^{\alpha}}}{2(A\rho+\frac{B\alpha}
{\rho^{\alpha}})}<1,
\end{equation}
we have $u_{\xi}>0$ from \eqref{eq6.6} for $A,B$ sufficiently small, so
the set $(u,\rho)$ which can be connected to $(u_-,\rho_-)$ by the
backward rarefaction wave is made up of the half-branch of
$\overleftarrow{R}(u_-,\rho_-)$ with $u\geq u_-$.

Similarly, for the forward  rarefaction wave, we have $u_{\xi}>0$
for $A,B$ sufficiently small, and the set $(u,\rho)$ which can be
joined to $(u_-,\rho_-)$ by the forward rarefaction wave is made up
of the half-branch of $\overrightarrow{R}(u_-,\rho_-)$ with $u\geq
u_-$.

Taking the limit $\rho\rightarrow 0$ in the second equation of \eqref{eq6.4}, it follows that
\begin{equation}\label{eq6.8}
\lim_{\rho\rightarrow
0}\sqrt{u}=\sqrt{u_-}+\frac{1}{2}\int^{\rho_-}_{0}
\frac{\sqrt{As+\frac{B\alpha}{s^{\alpha}}}}{s} ds.
\end{equation}
Since $\underset{s\rightarrow 0}\lim
\Big(s^{1+\frac{\alpha}{2}}\frac{\sqrt{As+\frac{B\alpha}{
s^{\alpha}}}}{s}\Big)=\sqrt{B \alpha}$, the integral $
\int^{\rho_-}_{0} \frac{\sqrt{As+\frac{B\alpha}{s^{\alpha}}}}{s} ds$ is
divergent owing to Cauchy criterion for integral of an unbounded
function. Thus, for the backward rarefaction wave, from \eqref{eq6.8} one can get
$\underset{\rho\rightarrow 0}\lim u=+\infty$ .

Performing the limit $\rho\rightarrow +\infty$ in the second
equation of \eqref{eq6.5} yields
\begin{equation}\label{eq6.9}
\lim_{\rho\rightarrow
+\infty}\sqrt{u}=\sqrt{u_-}+\frac{1}{2}\int^{+\infty}_{\rho_-}
\frac{\sqrt{As+\frac{B\alpha}{s^{\alpha}}}}{s} ds.
\end{equation}
Since
\begin{equation}\label{eq6.10}
\dis \int^{+\infty}_{\rho_-} \frac{\sqrt{As+\frac{B\alpha}{s^{\alpha}}}}{s}
ds>\dis \int^{+\infty}_{\rho_-} \frac{\sqrt{As}}{s}
ds=+\infty.
\end{equation}
Thus, for the forward rarefaction wave, we have
$\underset{\rho\rightarrow +\infty}\lim u=+\infty$ from \eqref{eq6.9}
according to comparison test for infinite integral.

For a bounded discontinuity at  $\xi=\overline{\sigma}^{AB}$,  the
Rankine-Hugoniot relation
\begin{align}\label{eq6.11}
\left\{\begin{array}{l}
-\overline{\sigma}^{AB}[\rho]+[\rho u]=0,\cr\noalign {\vskip3truemm}
-\overline{\sigma}^{AB}[\rho
u+\dis\frac{A}{2}\rho^2-\frac{B}{1-\alpha}\rho^{1-\alpha}]+[\rho
u^2+A\rho^2 u-B\rho^{1-\alpha}u]=0
\end{array}\right.
\end{align}
holds.

Eliminating $\overline{\sigma}^{AB}$ from \eqref{eq6.11}, we obtain
\begin{equation}\label{eq6.12}
(u_r-u_l)^2=E_1(u_l,\rho_l,u_r,\rho_r,A,B),
\end{equation}
where
\begin{align}\label{eq6.13}
\begin{array}{ll}
&E_1(u_l,\rho_l,u_r,\rho_r,A,B)=\dfrac{A}{2\rho_l}\rho_r^2
u_r+\dfrac{\alpha B}{(1-\alpha)\rho_l}\rho_r^{1-\alpha}
u_r+\dfrac{A}{2\rho_r}\rho_l^2 u_l+\dfrac{\alpha B}{(1-\alpha)
\rho_r}\rho_l^{1-\alpha} u_l
\\[5mm]
&-A\rho_r(u_r-\dfrac{1}{2}u_l)-\dfrac{B}{\rho_l^{\alpha}}(\dfrac{1}{(1-\alpha)}u_r-u_l)+A\rho_l(\dfrac{1}{2}u_r-u_l)+\dfrac{B}{\rho_r^{\alpha}}(u_r-\dfrac{1}{1-\alpha}u_l).
\end{array}
\end{align}
It is easy to check that $E_1(u_l,\rho_l,u_r,\rho_r,A,B)>0$.

Using the Lax entropy inequalities, one can get the backward shock
wave satisfies
\begin{equation}\label{eq6.14}
\overline{\sigma}^{AB}<\overline{\lambda_1}(u_l,\rho_l), \quad
\overline{\lambda_1}(u_r,\rho_r)<\overline{\sigma}^{AB}<\overline{\lambda_2}(u_r,\rho_r),
\end{equation}
and the forward shock wave satisfies
\begin{equation}\label{eq6.15}
\overline{\lambda_1}(u_l,\rho_l)<\overline{\sigma}^{AB}<\overline{\lambda_2}(u_l,\rho_l),
\quad
\overline{\lambda_2}(u_r,\rho_r)<\overline{\sigma}^{AB}.
\end{equation}

Then, from \eqref{eq6.14}, we can obtain that the following inequality holds
for backward shock wave
\begin{equation}\label{eq6.16}
\dis\frac{-\sqrt{u_r(A\rho_r+\frac{B\alpha}{\rho_r^{\alpha}})}}{\rho_l}<\frac{u_r-u_l}{\rho_r-\rho_l}<
\frac{-\sqrt{u_l(A\rho_l+\frac{B\alpha}{\rho_l^{\alpha}})}}{\rho_r}.
\end{equation}
Associating with $\lambda_1$, \eqref{eq6.16} implies that $\rho_l<\rho_r$
and $u_r<u_l$. In a analogous way, for the forward shock wave, we can deduce that
$\rho_l>\rho_r$ and $u_r<u_l$ from \eqref{eq6.15}.

Now, given a left state $(u_-,\rho_-)$, we have the following
backward shock wave curve
\begin{equation}\label{eq6.17}
 \overleftarrow{S}(u_-,\rho_-):\
u-u_-=-\sqrt{E_1(u_-,\rho_-,u,\rho,A,B)}, \ \ \ \rho>\rho_-,
\end{equation}
and the forward shock wave curve
\begin{equation}\label{eq6.18}
 \overrightarrow{S}(u_-,\rho_-):\
u-u_-=-\sqrt{E_1(u_-,\rho_-,u,\rho,A,B)},\ \ \
\rho<\rho_-.
\end{equation}

For the backward shock wave, through differentiating $u$ with
respect to $\rho$ in \eqref{eq6.17}, it follows that
\begin{equation}\label{eq6.19}
\bigg(1+\frac{E_2(u_-,\rho_-,u,\rho,A,B)}{2\sqrt{E_1(u_-,\rho_-,u,\rho,A,B)}}\bigg)u_{\rho}=\frac{-E_3(u_-,\rho_-,u,\rho,A,B)}{2\sqrt{E_1(u_-,\rho_-,u,\rho,A,B)}},
\end{equation}
where
\begin{equation}\label{eq6.20}
E_2(u_-,\rho_-,u,\rho,A,B)=\frac{A}{2\rho_-}\rho^2+\frac{\alpha
B}{(1-\alpha)\rho_-}\rho^{1-\alpha}-A\rho-\frac{B}{(1-\alpha)\rho_-^{\alpha}}+\frac{A}{2}\rho_-+\frac{B}{\rho^{\alpha}},
\end{equation}
\begin{align}\label{eq6.21}
\begin{array}{ll}
\hspace{-0.1cm}E_3(u_-,\rho_-,u,\rho,A,B)=Au
(\dfrac{\rho}{\rho_-}-1)+\dfrac{\alpha
Bu}{\rho^{\alpha}}(\dfrac{1}{\rho_-}-\dfrac{1}{\rho})+\dfrac{\alpha
B u_-}{(1-\alpha)\rho^{2}}(\rho^{1-\alpha}-\rho_-^{1-\alpha})\\[5mm]
\hspace{3.5cm}+\dfrac{A}{2}\big(1-(\dfrac{\rho_-}{\rho})^2\big)u_-.
\end{array}
\end{align}
Thus, for the backward shock wave, we have $u_{\rho}<0$  from \eqref{eq6.19}
as $A,B$ is sufficiently small. Similarly, for the forward shock
wave, we can obtain $u_{\rho}>0$ as $A,B$ is sufficiently small.

For the backward shock wave $\overleftarrow{S}(u_-,\rho_-)$,
taking $u=0$ in \eqref{eq6.17} gives
\begin{equation}\label{eq6.22}
u_-=-A\rho_-+\frac{B}{\rho_-^{\alpha}}+\frac{A}{2}\rho
-\frac{B}{(1-\alpha)\rho^{\alpha}}+\frac{A}{2\rho}\rho_-^2+\frac{\alpha
B}{(1-\alpha)\rho}\rho_-^{1-\alpha}.
\end{equation}

Let
\begin{equation}\label{eq6.23}
f(\rho)=-A\rho_-+\frac{B}{\rho_-^{\alpha}}+\frac{A}{2}\rho
-\frac{B}{(1-\alpha)\rho^{\alpha}}+\frac{A}{2\rho}\rho_-^2+\frac{\alpha
B}{(1-\alpha)\rho}\rho_-^{1-\alpha}-u_-.
\end{equation}
Then $f(\rho)$ is continuous with respect to $\rho$ and
$f(\rho_-)f(+\infty)<0$. Thus, there exist
$\rho_0\in[\rho_-,+\infty)$ such that $f(\rho_0)=0$, which implies
that the backward shock wave curve intersects with the
$\rho$-axis at a point. From \eqref{eq6.18}, it is not difficult to check
that $\overrightarrow{S}(u_-,\rho_-)$ has a intersection point with the
$\rho$-axis.

From the analysis above, fixing a left state $(u_-,\rho_-)$, the
phase plane ($u,\rho>0$) can be divided into four regions by
the wave curves (see Fig. 2), denoted by
$\overleftarrow{S}\overrightarrow{S}(u_-,\rho_-)$,
$\overleftarrow{S}\overrightarrow{R}(u_-,\rho_-)$,
$\overleftarrow{R}\overrightarrow{S}(u_-,\rho_-)$ and
$\overleftarrow{R}\overrightarrow{R}(u_-,\rho_-)$, respectively.

Now, according to the right state $(u_+,\rho_+)$ in the different
regions, we can get four kinds of configurations of solutions.
Particularly, when $(u_+,\rho_+)\in
\overleftarrow{S}\overrightarrow{S}(u_-,\rho_-)$, the Riemann
solution contains two shock waves and a intermediate
constant states whose density may become singular as
$A,B\rightarrow0$. When $(u_+,\rho_+)\in
\overleftarrow{R}\overrightarrow{R}(u_-,\rho_-)$, the Riemann
solution contains two rarefaction waves and a nonvacuum intermediate
constant states that may be a vacuum state as $A,B\rightarrow0$.
Since the other two regions
$\overleftarrow{S}\overrightarrow{R}(u_-,\rho_-)$ and
$\overleftarrow{R}\overrightarrow{S}(u_-,\rho_-)$ have empty
interiors when $A,B\rightarrow0$, it suffices to study the limit
process for the two cases $(u_+,\rho_+)\in
\overleftarrow{S}\overrightarrow{S}(u_-,\rho_-)$ and
$(u_+,\rho_+)\in \overleftarrow{R}\overrightarrow{R}(u_-,\rho_-)$.

\begin{center}
\includegraphics*[154,497][407,639]{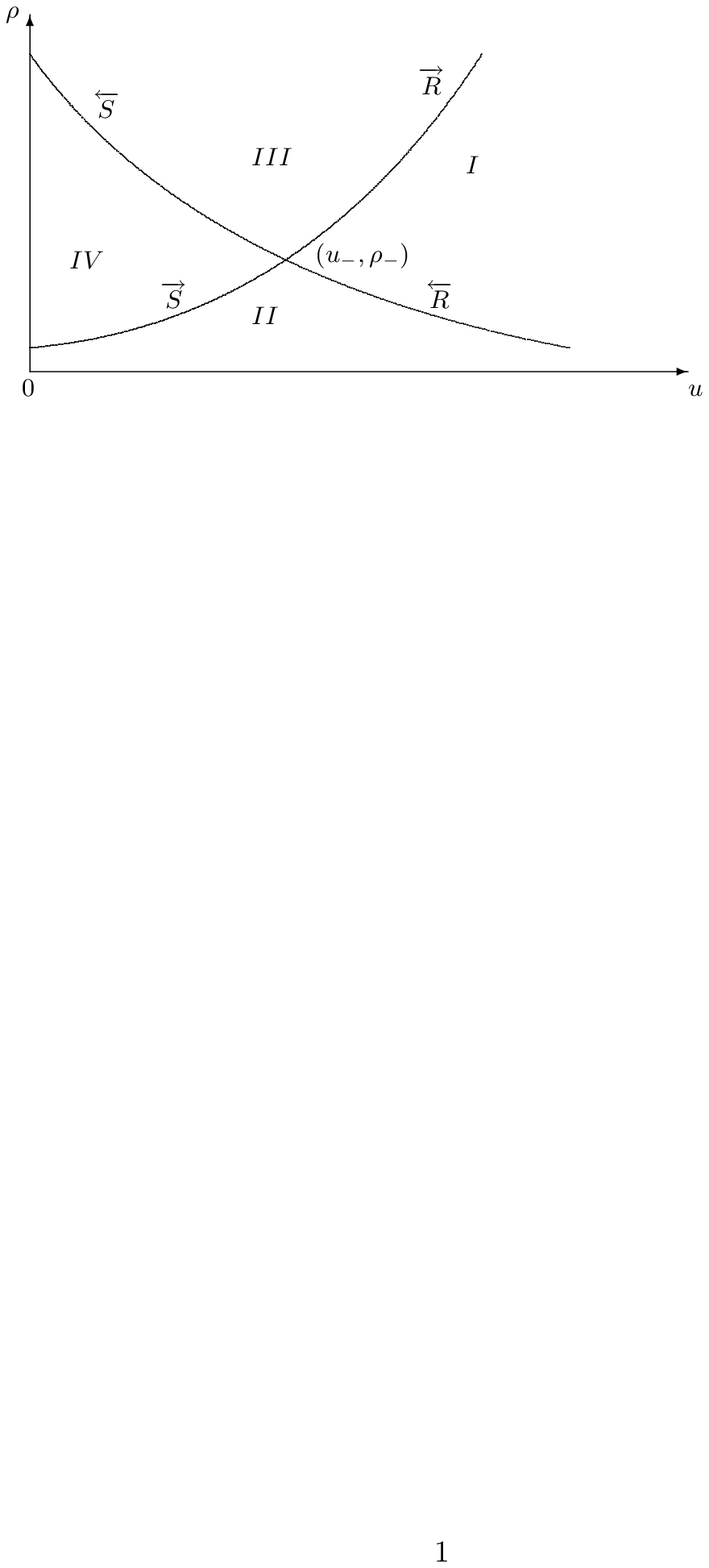}

{\text{ Fig. 2.} Curves of elementary waves.}

\end{center}

 \section{Formation of delta-shocks in solutions of \eqref{eq1.7} and \eqref{eq1.8}}

This section analyzes the limits as $A,B\rightarrow0$ of the
solutions of \eqref{eq1.7} and \eqref{eq1.8} in the case $(u_+,\rho_+)\in
\overleftarrow{S}\overrightarrow{S}(u_-,\rho_-)$ with $u_+<u_-$.

 \subsection{Limit behaviour of the Riemann solutions as $A,B\rightarrow0$}

As stated previously, for any fixed $A,B>0$, since  $(u_-,\rho_-)$ and
$(u_*^{AB},\rho_*^{AB})$  are connected by a backward shock wave
$\overleftarrow{S}$ with the propagation speed
$\overline{\sigma_1}^{AB}$, $(u_*^{AB},\rho_*^{AB})$ and
$(u_+,\rho_+)$ are connected by a forward shock wave
$\overrightarrow{S}$ with the propagation speed
$\overline{\sigma_2}^{AB}$, respectively. Thus we have
\begin{equation}\label{eq7.1}
\dis u_*^{AB}-u_-=-\sqrt{E_1(u_-,\rho_-,u_*^{AB},\rho_*^{AB},A,B)},\
\ \ \ \ \ \ \rho_-<\rho_*^{AB}
\end{equation}
on $\overleftarrow{S}$, and
\begin{equation}\label{eq7.2}
\dis u_+-u_*^{AB}=-\sqrt{E_1(u_*^{AB},\rho_*^{AB},u_+,\rho_+,A,B)},\
\ \ \ \ \ \ \rho_+<\rho_*^{AB}
\end{equation}
on $\overrightarrow{S}$.

Then we have the following lemmas.

\vspace{0.2cm}

 \noindent{\textbf{Lemma 7.1.}}\textit{
$\lim\limits_{A,B\rightarrow0}\rho_*^{AB}=+\infty$.}

\vspace{0.2cm}

The proof of this lemma is similar to that of Lemma 4.1. In fact, it
follows from \eqref{eq7.1} and \eqref{eq7.2} that
\begin{equation}\label{eq7.3}
\dis
u_+-u_-=-\Big(\sqrt{E_1(u_-,\rho_-,u_*^{AB},\rho_*^{AB},A,B)}+\sqrt{E_1(u_*^{AB},\rho_*^{AB},u_+,\rho_+,A,B)}\Big).
\end{equation}
Note that $u_+<u_*^{AB}<u_-$. We let $A,B\rightarrow0$ in \eqref{eq7.3} and
employ the method of proof by contradiction to obtain the desired
result.

\vspace{0.2cm}

\noindent {\textbf{Lemma 7.2.}} Set $\sigma=\dis
\frac{\sqrt{\rho_-}u_-+\sqrt{\rho_+}u_+}{\sqrt{\rho_-}+\sqrt{\rho_+}}$.
Then \textit{
\begin{equation}\label{eq7.4}
\lim\limits_{A,B\rightarrow0}u_*^{AB}=\lim\limits_{A,B\rightarrow0}\overline{\sigma_1}^{AB}=\lim\limits_{A,B\rightarrow0}\overline{\sigma_2}^{AB}=
\sigma,
\end{equation}
\begin{equation}\label{eq7.5}
\lim\limits_{A,B\rightarrow0}A\rho_*^{AB}=0.
\end{equation}}

\vspace{0.2cm}

\textbf{Proof.} Passing to the limit $A,B\rightarrow0$ in \eqref{eq7.3} and
noting Lemma 7.1, we reach
\begin{equation}\label{eq7.6}
\lim\limits_{A,B\rightarrow0}\sqrt{\frac{A}{2}(\rho_*^{AB})^2u_*^{AB}+\frac{\alpha
B}{1-\alpha}(\rho_*^{AB})^{1-\alpha}u_*^{AB}}=\frac{\sqrt{\rho_-\rho_+}(u_--u_+)}{\sqrt{\rho_-}+\sqrt{\rho_+}}.
\end{equation}
Taking $A,B\rightarrow0$ in \eqref{eq7.1}, one can obtain that
$\lim\limits_{A,B\rightarrow0}u_*^{AB}=\sigma$.

From \eqref{eq6.11}, $\overline{\sigma_1}^{AB}$ and
$\overline{\sigma_2}^{AB}$ can be calculated by
\begin{equation}\label{eq7.7}
\dis\overline{\sigma_1}^{AB}=\frac{\rho_*^{AB}u_*^{AB}-\rho_-u_-}{\rho_*^{AB}-\rho_-},\qquad
\dis\overline{\sigma_2}^{AB}=\frac{\rho_+u_+-\rho_*^{AB}u_*^{AB}}{\rho_+-\rho_*^{AB}},
\end{equation}
thus one can easily check that \eqref{eq7.4} holds. Then \eqref{eq7.4} and \eqref{eq7.6}
yield \eqref{eq7.5}.

\vspace{0.2cm}

Combining \eqref{eq7.7} with Lemmas 7.1-7.2, we can obtain the following
lemma.

\vspace{0.2cm}

\noindent {\textbf{Lemma 7.3.}} \textit{$
\lim\limits_{A,B\rightarrow0}\rho_*^{AB}(\overline{\sigma_2}^{AB}-\overline{\sigma_1}^{AB})=\sigma[\rho]-[\rho
u]. $}

\subsection{Weighted delta shock waves}

Now, we  show the theorem characterizing the limit as
$A,B\rightarrow0$ for the case $u_+< u_-$ and $(u_+,\rho_+)\in
\overleftarrow{S}\overrightarrow{S}(u_-,\rho_-)$.

\vspace{0.2cm}

\noindent {\textbf{Theorem 7.4.}} \textit{Let  $u_+< u_-$.  Assume
$(u^{AB},\rho^{AB})$ is a two-shock wave solution of \eqref{eq1.7} and
\eqref{eq1.8} constructed in Section 6. Then, when $A,B\rightarrow0$,
$\rho^{AB}$ and $\rho^{AB}u^{AB}$ converge in the sense of
distributions, and the limit functions of $\rho^{AB}$ and $\rho^{AB}
u^{AB}$ are the sums of a step function and a $\delta$-function with
the weights \
$$\dis\frac{t}{\sqrt{1+\sigma^2}}(\sigma[\rho]-[\rho u])\ \ and \ \ \dis\frac{t}{\sqrt{1+\sigma^2}}(\sigma[\rho u]-[\rho u^2]),$$
respectively, which form a delta shock solution of \eqref{eq1.4} with the
same Riemann data \eqref{eq1.8}.}

\vspace{0.2cm}

\textbf{Proof.} (i). Set $\xi=x/t$. Then for each $A,B>0$, the
 Riemann solution containing $\overleftarrow{S}$ and $\overrightarrow{S}$ can be expressed as
\begin{align}\label{eq7.8}
(u^{AB},\rho^{AB})(\xi)=\left\{\begin{array}{ll}
 (u_-,\rho_-),&\xi<\overline{\sigma_1}^{AB}, \\[2mm]
 (u_*^{AB},\rho_*^{AB}), &\overline{\sigma_1}^{AB}<\xi<\overline{\sigma_2}^{AB},\\[2mm]
 (u_+,\rho_+), &\xi>\overline{\sigma_2}^{AB},
 \end{array}\right.
\end{align}
satisfying weak formulations: For any $\phi\in
C^1_0(-\infty,+\infty)$,
\begin{equation}\label{eq7.9}
\dis\int^{+\infty}_{-\infty}\rho^{AB}( u^{AB}-
\xi)\phi'd\xi-\dis\int^{+\infty}_{-\infty}\rho^{AB}\phi
d\xi=0,
\end{equation}
and
\begin{align}\label{eq7.10}
\begin{array}{l}
\hspace{-1.0cm}-\dis{\int^{+\infty}_{-\infty}}\rho^{AB}\bigg(u^{AB}+\frac{A}{2}\rho^{AB}-\frac{B}{(1-\alpha)(\rho^{AB})^{\alpha}}\bigg)\xi
\phi'd\xi\\[5mm]
\hspace{2.0cm}+\dis{\int^{+\infty}_{-\infty}}\rho^{AB}\bigg(
(u^{AB})^2+A
\rho^{AB}u^{AB}-\frac{B}{(\rho^{AB})^{\alpha}}u^{AB}\bigg)\phi'd\xi\\[5mm]
\hspace{-1.3cm}=\dis\int^{+\infty}_{-\infty}\rho^{AB}\bigg(u^{AB}+\frac{A}{2}\rho^{AB}-\frac{B}{(1-\alpha)(\rho^{AB})^{\alpha}}\bigg)
\phi d\xi.
\end{array}
\end{align}

(ii). Now we turn to computing the limit of $\rho^{AB}u^{AB}$ and
$\rho^{AB}$ by using the two weak formulations. The first integral
on the left hand side of \eqref{eq7.10} can be decomposed into
\begin{equation}\label{eq7.11}
-\bigg(\dis\int^{\overline{\sigma_1}^{AB}}_{-\infty}+\dis\int_{\overline{\sigma_1}^{AB}}^{\overline{\sigma_2}^{AB}}
+\dis\int^{+\infty}_{\overline{\sigma_2}^{AB}}\bigg)\bigg(\rho^{AB}u^{AB}+\frac{A}{2}(\rho^{AB})^2-\frac{B}{(1-\alpha)}(\rho^{AB})^{1-\alpha}\bigg)\xi\phi'd\xi.
\end{equation}
The limit of the sum of the first and last terms of \eqref{eq7.11} is
\begin{align}\label{eq7.12}
\begin{array}{l}
-\lim\limits_{A,B\rightarrow0}\dis\int^{\overline{\sigma_1}^{AB}}_{-\infty}\bigg(\rho^{AB}u^{AB}+\frac{A}{2}(\rho^{AB})^2-\frac{B}{(1-\alpha)}(\rho^{AB})^{1-\alpha}\bigg)\xi\phi'd\xi\\[4mm]
\hspace{2.5cm}-\lim\limits_{A,B\rightarrow0}\dis\int_{\overline{\sigma_2}^{AB}}^{+\infty}\bigg(\rho^{AB}u^{AB}+\frac{A}{2}(\rho^{AB})^2-\frac{B}{(1-\alpha)}(\rho^{AB})^{1-\alpha}\bigg)\xi\phi'd\xi\\[4mm]
=\sigma[\rho u]\phi(\sigma)+\dis\int^{+\infty}_{-\infty}H(\xi-\sigma)\phi
 d\xi
\end{array}
\end{align}
with
$$
H(\xi-\sigma)=\left\{\begin{array}{ll}
 \rho_-u_-,&\xi<\sigma,\\[1mm]
 \rho_+u_+,&\xi>\sigma.
 \end{array}\right.
$$

For the limit of the second term  of \eqref{eq7.11}, applying Lemmas
7.1-7.3, we deduce that
\begin{align}\label{eq7.13}
\begin{array}{ll}
&-\lim\limits_{A,B\rightarrow0}\dis\int_{\overline{\sigma_1}^{AB}}^{\overline{\sigma_2}^{AB}}\bigg(\rho^{AB}u^{AB}+\frac{A}{2}(\rho^{AB})^2-\frac{B}{(1-\alpha)}(\rho^{AB})^{1-\alpha}\bigg)\xi\phi'd\xi\\[5mm]
&=-(\sigma[\rho]-[\rho u])\sigma^2\phi'(\sigma).
\end{array}
\end{align}

Then combining \eqref{eq7.12} and \eqref{eq7.13}, we have
\begin{align}\label{eq7.14}
\begin{array}{ll}
-\lim\limits_{A,B\rightarrow0}\dis{\int^{+\infty}_{-\infty}}\rho^{AB}\bigg(u^{AB}+\frac{A}{2}\rho^{AB}-\frac{B}{(1-\alpha)(\rho^{AB})^{\alpha}}\bigg)\xi
\phi'd\xi\\[5mm]
=\sigma[\rho u]\phi(\sigma)-(\sigma[\rho]-[\rho
u])\sigma^2\phi'(\sigma)+\dis\int^{+\infty}_{-\infty}H(\xi-\sigma)\phi
d\xi.
\end{array}
\end{align}

By computing the limit of the second integral on the left hand side
of \eqref{eq7.10}, we obtain that
\begin{align}\label{eq7.15}
\begin{array}{ll}
\lim\limits_{A,B\rightarrow0}\dis{\int^{+\infty}_{-\infty}}\rho^{AB}\bigg(
(u^{AB})^2+A
\rho^{AB}u^{AB}-\frac{B}{(\rho^{AB})^{\alpha}}u^{AB}\bigg)\phi'd\xi\\[5mm]
=\lim\limits_{A,B\rightarrow0}\bigg(\dis\int^{\overline{\sigma_1}^{AB}}_{-\infty}+\dis\int_{\overline{\sigma_1}^{AB}}^{\overline{\sigma_2}^{AB}}
+\dis\int^{+\infty}_{\overline{\sigma_2}^{AB}}\bigg)\rho^{AB}\bigg(
(u^{AB})^2+A
\rho^{AB}u^{AB}-\frac{B}{(\rho^{AB})^{\alpha}}u^{AB}\bigg)\phi'd\xi\\[5mm]
=-[\rho u^2]\phi(\sigma)+(\sigma[\rho]-[\rho
u])\sigma^2\phi'(\sigma).
\end{array}
\end{align}

Substituting \eqref{eq7.14} and \eqref{eq7.15} into \eqref{eq7.10} yields
\begin{align}\label{eq7.16}
\begin{array}{l}
\lim\limits_{A,B\rightarrow0}\dis\int^{+\infty}_{-\infty}\rho^{AB}\bigg(u^{AB}+\frac{A}{2}\rho^{AB}-\frac{B}{(1-\alpha)(\rho^{AB})^{\alpha}}\bigg)
\phi d\xi=
\lim\limits_{A,B\rightarrow0}\dis\int^{+\infty}_{-\infty}\rho^{AB}
u^{AB}\phi d\xi \\
=(\sigma[\rho u]-[\rho
u^2])\phi(\sigma)+\dis\int^{+\infty}_{-\infty}H(\xi-\sigma)\phi
d\xi.
\end{array}
\end{align}

Similarly, from \eqref{eq7.9}, one can get that
\begin{equation}\label{eq7.17}
 \lim\limits_{A,B\rightarrow0}\dis\int^{+\infty}_{-\infty}\rho^{AB}\phi d\xi
 =(\sigma[\rho]-[\rho u])\phi(\sigma)+\dis\int^{+\infty}_{-\infty}\widetilde{H}(\xi-\sigma)\phi
 d\xi,
\end{equation}
where
$$
\widetilde{H}(\xi-\sigma)=\left\{\begin{array}{ll}
 \rho_-,&\xi<\sigma,\\[1mm]
 \rho_+,&\xi>\sigma.
 \end{array}\right.
$$

(iii). Finally, we analyze the limit of $\rho^{AB}u^{AB}$ and
$\rho^{AB}$ by tracking the time-dependence of the weights of the
$\delta$-measures as $A,B\rightarrow0$.

Taking \eqref{eq7.16} into account, we have, for any $\psi\in
C^\infty_0(R^+\times R)$,
\begin{align}\label{eq7.18}
\begin{array}{l}
\lim\limits_{A,B\rightarrow0}\dis\int^{+\infty}_{0}\dis\int^{+\infty}_{-\infty}\rho^{AB}(x/t)u^{AB}(x/t)\psi(t,x)dxdt\\[5mm]
 =\lim\limits_{A,B\rightarrow0}\dis\int^{+\infty}_{0}\dis\int^{+\infty}_{-\infty}\rho^{AB}(\xi)u^{AB}(\xi)\psi(t,\xi t)d(\xi t)dt\\[5mm]
 =\dis\int^{+\infty}_{0}(\sigma[\rho u]-[\rho
u^2])t\psi(t,\sigma
t)dt+\dis\int^{+\infty}_{0}\dis\int^{+\infty}_{-\infty}H(x-\sigma
 t)\psi(t,x)dxdt,
\end{array}
\end{align}
in which, by the definition \eqref{eq2.2}, we get
$$
 \dis\int^{+\infty}_{0}(\sigma[\rho u]-[\rho
u^2])t\psi(t,\sigma
t)dt=\Big<w_1(\cdot)\delta_S,\psi(\cdot,\cdot)\Big>
$$
with
$$
w_1(t)=\dis\frac{t}{\sqrt{1+\sigma^2}}(\sigma[\rho u]-[\rho u^2]).
$$

Similarly, we can show that
\begin{align}\label{eq7.19}
\begin{array}{l}
\lim\limits_{A,B\rightarrow0}\dis\int^{+\infty}_{0}\dis\int^{+\infty}_{-\infty}\rho^{AB}(x/t)\psi(t,x)dxdt\\[5mm]
 =\Big<w_2(\cdot)\delta_S,\psi(\cdot,\cdot)\Big>+\dis\int^{+\infty}_{0}\dis\int^{+\infty}_{-\infty}\widetilde{H}(x-\sigma
 t)\psi(t,x)dxdt
\end{array}
\end{align}
with
$$
w_2(t)=\dis\frac{t}{\sqrt{1+\sigma^2}}(\sigma[\rho]-[\rho u]).
$$

The proof of Theorem 7.4 is completed.

\section{Formation of vacuums in solutions of \eqref{eq1.7} and \eqref{eq1.8}}

This section studies the limits of the solutions of \eqref{eq1.7} and \eqref{eq1.8} in the
case $(u_+,\rho_+)\in
\overleftarrow{R}\overrightarrow{R}(u_-,\rho_-)$ with $u_+>u_-$ as
$A,B\rightarrow0$.

\vspace{0.2cm}

In this case, the solution consists of two rarefaction waves
$\overleftarrow{R}$, $\overrightarrow{R}$ and an intermediate state
$(u_*^{AB},\rho_*^{AB})$, besides two constant states
$(u_{\pm},\rho_{\pm})$. They satisfy
\begin{align}\label{eq8.1}
\overleftarrow{R}:\ \left\{\begin{array}{l}
 \xi=\overline{\lambda_1}=u-\sqrt{u(A\rho+\dfrac{B\alpha}{\rho^{\alpha}})},\cr\noalign {\vskip2truemm}
 \sqrt{u}-\sqrt{u_-}=-\dis\frac{1}{2}\int^{\rho}_{\rho_-} \frac{\sqrt{As+\frac{B\alpha}{s^{\alpha}}}}{s} ds,\
\ \ \ \ \ \ \rho_*^{AB}\leq\rho\leq\rho_-,
\end{array}\right.
\end{align}
and
\begin{align}\label{eq8.2}
\overrightarrow{R}:\ \left\{\begin{array}{l}
\xi=\overline{\lambda_2}=u+\sqrt{u(A\rho+\dfrac{B\alpha}{\rho^{\alpha}})},\cr\noalign
{\vskip2truemm}
 \sqrt{u_+}-\sqrt{u}=\dis\frac{1}{2}\int^{\rho_+}_{\rho} \frac{\sqrt{As+\frac{B\alpha}{s^{\alpha}}}}{s} ds,\
\ \ \ \ \ \ \rho_*^{AB}\leq\rho\leq\rho_+.
\end{array}\right.
\end{align}

We now conclude a theorem as follows.

\vspace{0.2cm}

\noindent {\textbf{Theorem 8.1.}} \textit{When $u_+>u_-$ and
$(u_+,\rho_+)\in \overleftarrow{R}\overrightarrow{R}(u_-,\rho_-)$,
the vacuum state occurs as $A,B\rightarrow0$, and two rarefaction
waves become two contact discontinuities connecting the constant
states $(u_\pm,\rho_\pm)$ and the vacuum state $(\rho=0)$. }

\vspace{0.2cm}

The proof of Theorem 8.1 is similar to that of Theorem 5.1, here we omit it.

\vspace{0.2cm}

In summary, the limit solution in this case can be expressed as
\eqref{eq2.1}, which is a solution to the transport equations containing two
contact discontinuities $\xi=x/t=u_\pm$ and a vacuum state in
between.

\vspace{0.2cm}
\noindent {\textbf{Remark.}} The processes of formation of delta shock waves and vacuum states
can be examined with some numerical results as $A$ and $B$
decrease. The numerical simulations will be added to our submission soon.

\end{document}